\documentclass[12pt,reqno]{amsart}
\usepackage{amsmath,amssymb}
\usepackage[dvips]{graphicx,color}
\usepackage{latexsym,amssymb}
\usepackage{mathrsfs}
\usepackage[abbrev]{amsrefs}

\newtheorem{thm}{Theorem}[section]

\newtheorem{prop}[thm]{Proposition}
\newtheorem{lem}[thm]{Lemma}

\numberwithin{equation}{section}
\numberwithin{thm}{section}

\begin{document}

\begin{center}\large \bf 
On analyticity up to the boundary for critical 
\\
quasi-geostrophic equation 
\end{center}

\footnote[0]
{


E-mail: t-iwabuchi@tohoku.ac.jp

}
\vskip5mm

\begin{center}
Tsukasa Iwabuchi 

\vskip2mm

Mathematical Institute, 
Tohoku University\\
Sendai 980-8578 Japan

\end{center}

\vskip5mm

\begin{center}
\begin{minipage}{135mm}
\footnotesize
{\sc Abstract. } 
We study the Cauchy problem for the surface quasi-geostrophic equations with the critical dissipation 
in the two dimensional half space under the homogeneous Dirichlet boundary condition. 
We show the global existence, the uniqueness and the analyticity of solutions, 
and the real analyticity up to the boundary is obtained. 
We will show a natural way to estimate the nonlinear term for functions satisfying the Dirichlet boundary condition.

\end{minipage}
\end{center}

\section{Introduction}

We consider the critical surface quasi-geostrophic equations in the half space. 
\begin{equation}\label{QG1}
\displaystyle 
 \partial_t \theta 
  + (u \cdot \nabla ) \theta + \Lambda _D \theta  =0, 
  \quad u= \nabla ^{\perp} \Lambda _D ^{-1} \theta, 
 \qquad  t > 0 , x \in \mathbb R ^2_+, 
\end{equation}
\begin{equation}\label{QG2}
 \theta(0,x) = \theta_0(x) , 
\qquad  x \in \mathbb R ^2_+, 
\end{equation}
where 
$\mathbb R^2_+ := \{ (x_1, x_2 ) \in \mathbb R^2 \, | \, x_n > 0 \}$, 
$\nabla^\perp := (-\partial_{x_2} , \partial_{x_1})$, 
$\Lambda_D$ is the square root of the Dirichlet Laplacian. 
The equations are known as an important model in geophysical fluid dynamics, 
which is derived from general quasi-geostrophic equations in the special case of 
constant potential vorticity and buoyancy frequency (see~\cite{La_1959,Pe_1979}). 
The purpose of this paper is to show the existence of global solutions 
for initial data in scaling critical Besov spaces, and the analyticity.

Let us recall 
existing results, where the space is the whole space. 
If we consider the fractional Laplacian of the order $\alpha$, 
$(-\partial_x^2) ^{\alpha/2}$, with $0 < \alpha \leq 2$, instead of 
the square root of the Laplacian, 
then 
the case when $\alpha < 1, \alpha = 1, \alpha > 1$ are called 
sub-critical case, critical case, super-critical case, respectively. 
It is known that the global-in-time regularity is obtained for the 
sub-critical case and the critical case. 
The sub-critical case can be treated, by $L^\infty$-maximum principle, 
and the critical case is delicate. In the critical case, 
the regularity with small data was proved by 
Constantin, Cordoba and Wu~\cite{CoCoWu-2001} 
(see also Constantin and Wu~\cite{ConWu-1999}).  
The poroblem for large data case was solved by Caffarelli and Vasseur~\cite{CaVa-2010}, 
Kiselev, Nazarov and Volberg~\cite{KNV-2007}. 
As another approach, Constantin and Vicol~\cite{CV-2012} proved the 
global regularity by 
nonlinear maximum principles 
in the form of a nonlinear lower bound on the fractional Laplacian. 
On the other hand, in the super-critical case, 
the regularity only for small data is known (see~\cite{CotVic-2016}), and 
blow-up for smooth solutions is an open problem.

In bounded domains with smooth boundary, the equations was introduced by 
Constantin and Ignatova~(\cite{CoIg-2016,CoIg-2017}). 
Local existence was shown by Constantin and Nguyen~\cite{CoNg-2018-2}, 
and global existence of weak solutions was proved 
by Constantin and Ignatova~\cite{CoIg-2017} 
(see also the paper by Constantin and Nguyen~\cite{CoNg-2018} for the inviscid case). 
An interesting question here is how to understand the behavior of the solutions. 
A priori bounds of smooth solutions was obtained by Constantin and 
Ignatova~\cite{CoIg-2016}, 
and interior Lipschitz continuity of weak solutions was studied by Ignatova~\cite{Ig-2019}. 
Recently, Constantin and Ignatova~\cite{CoIg-2020} considered 
the quotient of the solution by the first eigen function to 
investigate near the boundary, and gave a condition to obtain the global regularity 
up to the boundary. 
Stokols and Vasseur~\cite{StVa-2020} constructed global-in-time weak solutions 
with H\"older regularity up to the boundary. 
We should note from the viewpoint of smooth solutions 
that regularity near the boundary is guaranteed for a short time to the best of our knowlegde. 
As for the half space case, the possibility is pointed out in \cite{CoIg-2016} when 
the support of the initial data is away from the boundary.

In this paper, we consider the problem in the half space to show global-in-time 
regularity, and furthermore, analyticity in spacetime. 
Initial data in our theorem can have its support around the boundary. 
We will utilize the odd extention with respect to $x_2$ technically, 
but the reason of the half space is just for the sake of the simplicity. 
Our method is based on the one for the whole space, introducing 
Besov spaces associated with the Dirichlet Laplacian. 
Related idea to handle product estimate 
can be found in the paper~\cite{Iw:preprint}, where 
the validity of bilinear estimates for functions with the Dirichlet boundary 
condition is discussed. 
We should also remark that our domain, the half space, is one of the simplest domains, 
and the odd reflection and 
the existing result in the whole space $\mathbb R^2$ imply the existence of solutions below formally, 
but the main subject here is the behavior of functions on the boundary. 
We suppose to establish a method applicable to more general domains in the future, and 
the purpose is to state theorems in intrinsic framework.

We state two theorems. The first theorem concerns with the integral equation 
with the small data, where it seems easier to explain 
the proof near boundary clearly. 
The second theorem studies the data belonging to the largest scaling critical Besov space, where 
the smooth functions exist dense. 
We here introduce a formal definition of Besov spaces, 
which is defined precisely in section 2. 
For every function $f \in L^1 (\mathbb R^2_+ ) + L^\infty (\mathbb R^2_+)$, 
let $f_{odd}$ be defined by 
\[
f_{odd} (x_1,x_2) 
:= \begin{cases}
f(x_1,x_2) & \text{ if } x_2 > 0, 
\\
-f(x_1,-x_2) & \text{ if } x_2 < 0, 
\end{cases}
\]
where this extention is justfied as a locally integrable function at least. 
Norm of our Besov spaces can be understood 
through the spaces on the whole space by 
\[
\| f \|_{B^s_{p,q}(\Lambda_D)} = 
\| f_{odd} \|_{B^s_{p,q} (\mathbb R^2)}, \quad 
\| f \|_{\dot B^s_{p,q}(\Lambda_D)} = 
\| f_{odd} \|_{\dot B^s_{p,q} (\mathbb R^2)}. 
\]
We then have the following theorem.

\begin{thm}{\rm (}Solutions of the integral equation with small data{\rm)}
\label{thm:small}
Let $\theta _0 \in \dot B^0_{\infty,1}(\Lambda _D ) $ 
 be sufficiently small. 
 Then the integral equation
\[
\theta (t) = e^{- t \Lambda_D} \theta _0 
  -\int_0^t e^{-(t-\tau)\Lambda_D} 
   \Big( (u \cdot \nabla ) \theta \Big) ~d\tau , \quad 
   u= \nabla ^{\perp} \Lambda _D ^{-1} \theta
\]
posseses a unique global solution $\theta$ such that 
\begin{gather*}
\theta \in C([0,\infty ), \dot B^0_{\infty,1} (\Lambda _D)) 
\cap L^1 (0,\infty ; \dot B^1_{\infty,1} (\Lambda _D )). 
\end{gather*}
Furthoermore, there exists $C >0$ such that for any $\alpha , \beta_1, \beta_2 
\in \mathbb N\cup\{ 0 \}$, 
\[
t^{\alpha + \beta_1 + \beta _1 } \| \partial _t^\alpha \partial_{x_1} ^{\beta_1} \partial_{x_2} ^{\beta_2} \theta(t)  \| _{L^\infty (\mathbb R^2_+)}
\leq C^{\alpha + \beta_1 + \beta_2} \alpha ! \beta_1 ! \beta_2 ! , 
\]
and in particular, $\theta(t)$ is real analytic in space and time if $t > 0$. 
\end{thm}

%
\bigskip 

Next theorem concerns with the initial data in the space corresponding to 
the largest scaling critical Besov space $\dot B^0_{\infty,\infty}(\Lambda _D)$. 
We recall the the paper~\cite{WZ-2011} by Wang-Zhang, 
who take the data in the space defined by the completion of $C_0^\infty (\mathbb R^2)$. 
We can generalize it to the Besov space $B^0_{\infty,\infty} (\mathbb R^2)$ with  
taking the completion with small high frequency (see~\cite{Iw-2020}). 
Under this motivation, we have:

\begin{thm}\label{thm:1}
{\rm (}Solutions in the largest scaling critical Besov space{\rm )}
Let $\theta_0 \in B^0_{\infty,\infty} (\Lambda _D)$ be such that 
\begin{equation}\label{317-1}
\lim _{j \to \infty} \| \phi_j (\Lambda _D) \theta_0 \|_{L^\infty} = 0 , \quad 
\Big( 1 - \sum_{j \geq 0} \phi_j (\Lambda _D) \Big) \theta _0 \in \dot B^0_{\infty,1} (\Lambda_D) . 
\end{equation}
Then the problem \eqref{QG1} and \eqref{QG2} possess 
a unique solution $\theta$ such that 
\[
\theta \in C([0,\infty), B^0_{\infty,\infty} (\Lambda_D)) \cap 
L^1_{loc} (0,\infty ; \dot B^1_{\infty,\infty} (\Lambda_D)) , 
\]
\[
\Big( 1 - \sum_{j \geq 0} \phi_j (\Lambda _D) \Big)\theta \in C([0,\infty), \dot B^0_{\infty,1} (\Lambda_D)) .
\]
Furthoermore, $\theta(t)$ is real analytic in space and time if $t > 0$. 
\end{thm}

\noindent 
{\bf Remark}. 
(i) 
If we replace the first condition of \eqref{317-1} with the smallness of the high spectral 
component 
\[
\limsup_{j \to \infty} \| \phi_j(\lambda_D )\theta_0 \|_{L^\infty} \leq \delta ,  
\qquad \delta \ll 1,
\]
then we can also construct a unique solution such that 
$\theta$ is weak-* continuous in $B^0_{\infty,\infty}(\Lambda _D)$ with respect to $t \geq 0$. 
We can regard the smallness as the possibility of only small shock discussed in~\cite{CV-2012}.  

\noindent 
(ii) 
We impose the second condition of \eqref{317-1}  to justify $u = \nabla ^{\perp} \Lambda _D ^{-1} \theta 
\in \mathcal X'_D$. 

\noindent 
(iii) We will give a proof outline of global regularity based on the nonlinear muximum principle. 
Uniform bound in the H\"older spaces similar to Theorem~3.1 in \cite{Con-2017} enables 
us to repeat the fixed point argument in a time interval of fixed length.  


\vskip3mm

Let us give remarks for the proof of theorems. 
We apply a simple fixed point argument to the the proof of Theorem~\ref{thm:small}, 
as in \cite{Iw-2015} and \cite{Iw-2020}. 
To this end, we will derive bilinear estimates for $(u \cdot \nabla) \theta$, 
which is crucial to understand how to estimate near the boundary 
and contains the main idea of this paper.  
As for Theorem~\ref{thm:1}, we explain only proof outline, since the main idea 
near the boundary is same as in the proof of the first theorem and 
we can apply the proof in~\cite{Iw-2020}. 
We also mention the linear estimates which was established for more general 
framework (see~ \cite{Iw-2018}). 

We here focus on the discussion of the validity of the bilinear estimate for 
$(u \cdot \nabla) \theta$, 
how to obtain smoothness measured by the Dirichlet Laplacian up to the boundary. 
A standard argument would be to show 
\[
\| (u \cdot \nabla) \theta \|_{\dot H^s_p} 
\leq C ( \| u\|_{\dot H^s_{p_1}} \| \nabla \theta \|_{p_2}  
   + \| u \|_{L^{p_3}} \| \nabla\theta \|_{\dot H^s_{p_4}}) , 
   \quad s>0, \frac{1}{p} = \frac{1}{p_1} + \frac{1}{p_2} = \frac{1}{p_3} + \frac{1}{p_4},
\]
and apply the bounddedness of the Riesz transform. However,  this method causes 
a problem from the boundary value of the functions in the right hand side 
because of $x_2$ derivative, which yields 
$\| \partial _{x_2}\Lambda_D^{-1} \theta \|_{\dot H^s_{p_1}}, 
\| \partial _{x_2} \theta \|_{\dot H^s_{p_4}} = \infty$ for large $s$, in general. 
On the other hand, 
we investigate the boundary value of $(u \cdot \nabla) \theta$ itself 
which leads to a natural estimate. We can then have a estimate in Besov spaces 
(see~Propositiomn~\ref{prop:0309}). 
\[
\| (u \cdot \nabla) \theta \|_{\dot B^s_{p,1}(\Lambda_D)} 
\leq C ( \| \theta \|_{\dot B^s_{p_1,1}(\Lambda_D)} \| \theta \|_{\dot B^1_{p_2,1}(\Lambda_D)  }
   + \| \theta \|_{\dot B^0_{p_3,1}(\Lambda_D)} \| \theta \|_{\dot B^s_{p_4,1}(\Lambda_D)}) . 
\]
In a word, the most important point is: 
\textit{If $\theta$ satisfies the Dirichlet boundary condition, then 
$(u \cdot \nabla) \theta$ does.}
It would be possible to apply this argument not only to the half space  
but also to more general domains with smooth boundary. 
One can also find how derivatives affect boundary value of functions with 
the Dirichlet and the Neumann Laplacian in the papers~\cite{Iw-2018-2,Iw:preprint}, 
where the validity of product estimate for $fg$ is discussed and we can not expect it 
when $s > 2 + 1/p$. Briefly speaking, 
$\partial_{x_1}$ maps from $\dot B^{1}_{p,q} (\Lambda _D)$ to $\dot B^0_{p,q}(\Lambda_D)$, 
but $\partial _{x_2} $ maps to $\dot B^1_{p,q}(\Lambda_D)$ to $\dot B^0_{p,q}(\Lambda _N)$, 
where $\Lambda _N$ is the square root of the Neumann Laplacian.

\bigskip

This paper is organized as follows. 
In section 2, we recall the definition of Besov spaces associated with the Dirichlet 
Laplacian, and several properties to study the boundary value of functions. 
We prove Theorem~\ref{thm:small} in section 3, and explain idea of proof of Theorem~\ref{thm:1} 
in section 4.

\bigskip 

\noindent 
{\bf Notations.} We denote by $-\Delta _D$ the Dirichlet Laplacian on $L^2 (\Omega)$.  
We write $x=(x_1,x_2)$. 
For any function $f$ on $\mathbb R^2_+$, $f_{odd}, f_{even}$ are odd and even extentions of $f$ 
with respect to $x_2$, 
\[
f_{odd} (x) 
= \begin{cases}
f(x) & \text{ if } x_2 > 0, 
\\
-f(x_1,-x_2) & \text{ if } x_2 < 0, 
\end{cases}
\quad 
f_{even} (x) 
= \begin{cases}
f(x) & \text{ if } x_2 > 0, 
\\
f(x_1,-x_2) & \text{ if } x_2 < 0. 
\end{cases}
\]
For multi-index $\alpha= (\alpha _1 , \alpha _2) \in (\mathbb N \cup \{ 0 \})^2$,  
let 
$\partial_x^\alpha = \partial_{x_1}^{\alpha_1} \partial _{x_2}^{\alpha_2}$ and 
$|\alpha| = \alpha_1 + \alpha _2$. 
%
We denote by $\Delta_{\mathbb R^2}$ the Laplacian on $\mathbb R^2$, 
$G_t = G_t (x)$ the Gauss kernel 
$G_t (x) = (4\pi t)^{-1} e^{-\frac{|x|^2}{4t}}$, and 
$e^{t\Delta _D}$ the semigroup generated by the Dirichlet Laplacian 
\[
e^{t\Delta_D} f = \int_{\mathbb R^2_+} 
\big( G_t (x-y) - G_t (x+y) \Big) f(y) ~dy . 
\]
Let $\{ \phi_j \}_{j \in \mathbb Z}$ be the dyadic decomposition of the 
unity such that $\phi_j$ is a non-negative function in $C_0 ^\infty (\mathbb R)$ and 
\[
{\rm supp \, } \phi_0 \subset [2^{-1} , 2], \quad 
\phi_j (\lambda) = \phi_0 \left( \dfrac{\lambda}{2^j} \right) , 
\quad 
\sum _{j \in \mathbb Z} \phi_j (\lambda)  = 1  
\text{ for any } \lambda > 0 . 
\]
Let $\psi \in C_0^\infty (\mathbb R)$  be a non-negative function such that 
\[
{\rm supp \, } \psi \subset (-\infty , 2], \quad 
\psi(\lambda) + \sum_{j=1}^\infty \phi_j(\lambda) = 1 \text{ for any } \lambda \geq 0.
\]
We write the Fourier transform and the inverse Fourier transform. 
\[
\mathcal F[f](\xi) 
= \frac{1}{2\pi} \int_{\mathbb R^2} e^{ - i \xi \cdot x} f(x) dx, 
\qquad 
\mathcal F^{-1}[f](x) 
= \frac{1}{2\pi} \int_{\mathbb R^2} e^{ix \cdot \xi} f(\xi) d\xi . 
\]
The convolution $f*g$ is defined by the standard integral on $\mathbb R^2$. 
\[
f*g(x) = \int_{\mathbb R^2} f(x-y) g(y) dy.
\]
We use the following notations for norms of spaces in space and time as follows. 
\[
\begin{split}
\| f \|_{B^s_{p,q}(\Lambda_D)} = 
&
\| \psi(\Lambda_D) f \|_{L^p (\mathbb R^2_+)} 
+ 
\Big\|  
 \Big\{ 2^{sj} \| \phi_j(\Lambda _D) f \|_{L^p(\mathbb R^2_+)} 
 \Big\}_{j \in \mathbb N}
\Big\| _{\ell^q(\mathbb N)},
\\
\| f \|_{\dot B^s_{p,q}(\Lambda_D)} = 
&
\Big\|  
 \Big\{ 2^{sj} \| \phi_j(\Lambda _D) f \|_{L^p(\mathbb R^2_+)} 
 \Big\}_{j \in \mathbb Z}
\Big\| _{\ell^q(\mathbb Z)},
\\
\| f \|_{L^r(0,\infty; X)} 
= 
&\big\| \| f(t) \|_{X} \big\|_{L^r(0,\infty)} , 
\qquad X = L^p (\mathbb R^2_+), B^s_{p,q} (\Lambda _D ) , \dot B^s_{p,q} (\Lambda _D), 
\\
\| f \|_{\tilde L^r(0,\infty; \dot B^s_{p,q}(\Lambda _D))} 
= 
& \Big\|  
 \Big\{ 2^{sj} \| \phi_j(\Lambda _D) f \|_{L^r (0,\infty;L^p(\mathbb R^2_+))} 
 \Big\}_{j \in \mathbb Z}
\Big\| _{\ell^q(\mathbb Z)}.
\end{split}
\]

\section{Preliminary}

In this section, we recall the definition of the Besov spaces (see~\cite{IMT-2019}),  
the boundedness of the spectral multipliers (see~\cite{IMT-2018,Ou_2005,ThOuSi-2002}), 
several lemmas to justify an argument by 
the odd or the even extention of $f$ and its derivatives. 

We start by defining the Dirichlet Laplacian $-\Delta _D$, and 
spaces of test functions of non-homogeneous type, $\mathcal X_D$, 
and of homogeneous type $\mathcal Z_D$, 
and their duals.

\bigskip

\noindent 
{\bf Definition}. 
(i)
Let $-\Delta _D $ be the Dirichlet Laplacian on $L^2 (\mathbb R^2_+)$ defined by 
\[
\begin{cases}
D(-\Delta _D) := \{ f \in H^1_0(\mathbb R^2_+)  \, | \, 
 \Delta f \in L^2 (\mathbb R^2_+)\}, 
\\
-\Delta _D f := -\Delta f =
- \left( \dfrac{\partial^2}{\partial x_1 ^2} f 
  + \dfrac{\partial^2}{\partial x_2 ^2} f 
  \right) , 
\quad f \in D(-\Delta _D). 
\end{cases}
\]

\noindent 
(ii) Let $\mathcal X_D = \mathcal X_D(\mathbb R^2_+)$ 
be a space of test functions of non-homogeneous type such that 
\[
\mathcal X_D := \{ f \in L^1(\mathbb R^2_+) \cap L^2 (\mathbb R^2_+)  
   \, | \, 
   p_m (f) < \infty \text{ for all } m \in \mathbb N 
\}, 
\]
where 
\[
p_m (f) := \| f \|_{L^1} 
  + \sup_{ j \in \mathbb N} 2^{mj} \| \phi_j (\Lambda _D) f \|_{L^1} . 
\]

\noindent 
(iii) Let $\mathcal Z _D = \mathcal Z_D (\mathbb R^2_+) $ 
be a space of test functions of homogeneous type such that 
\[
\mathcal Z_D := \{ f \in \mathcal X_D \, | \, 
   q_m (f) < \infty \text{ for all } m \in \mathbb N  \}, 
\]
where 
\[
q_m(f):= 
p_m (f) + \sup _{ j \leq 0} 2^{m |j|} \| \phi_j (\Lambda _D) f \|_{L^1} . 
\]

\noindent 
(iv) Let $\mathcal X_D'$, $\mathcal Z_D '$ be the topological duals of 
$\mathcal X_D , \mathcal Z_D$, respectively.

\bigskip 

It was proved in \cite{IMT-2019} that 
the spaces $\mathcal X_D $, $\mathcal Z_D$ are Fr\'echet spaces, 
and can regard their duals $\mathcal X_D', \mathcal Z_D'$  
as distribution spaces of non-homogeneous type and homogeneous type, 
respectively, 
which are variants of the space of the tempered distributions and the quatient space by the polynomials. 
We define Besov spaces associated with the Dirichlet Laplacian 
on the half space as follows. 

\bigskip

\noindent 
{\bf Definition}. 
Let $s \in \mathbb R$ and $1 \leq p,q \leq \infty$. 
\begin{enumerate}
\item[(i)] (Non-homogeneous Besov space) $B^s_{p,q}(\Lambda_D)$  is defined by 
\[
B^s_{p,q}(\Lambda_D) 
= \{ f \in \mathcal X_D ' \, | \, 
   \| f \|_{B^s_{p,q}(\Lambda_D)} < \infty \},
\]
where 
\[
\| f \|_{B^s_{p,q}(\Lambda_D)} := 
\| \psi(\Lambda_D) f \|_{L^p (\mathbb R^2_+)} 
+ 
\Big\|  
 \Big\{ 2^{sj} \| \phi_j(\Lambda _D) f \|_{L^p(\mathbb R^2_+)} 
 \Big\}_{j \in \mathbb N}
\Big\| _{\ell^q(\mathbb N)}.
\]

\item[(ii)] (Homogenous Besov space) 
$\dot B^s_{p,q}(\Lambda_D)$ is defined by 
\[
\dot B^s_{p,q}(\Lambda_D) :=  
\{ f \in \mathcal Z_D' \, | \, 
   \| f \|_{\dot B^s_{p,q}(\Lambda_D) } < \infty\}, 
\]
where 
\[
\| f \|_{\dot B^s_{p,q}(\Lambda_D)} 
:= 
\Big\|  
 \Big\{ 2^{sj} \| \phi_j(\Lambda _D) f \|_{L^p(\mathbb R^2_+)} 
 \Big\}_{j \in \mathbb Z}
\Big\| _{\ell^q(\mathbb Z)}.
\]
\end{enumerate}

\bigskip 

It is proved in \cite{IMT-2019} that the Besov spaces $B^s_{p,q}(\Lambda _D)$, $\dot B^s_{p,q}(\Lambda _D) $ 
are Banach spaces and satisfy standard properties such as  lift properties, 
embedding theorems of Sobolev type 
as well as the whole space case (see~\cite{Tri_1983}). 
We here recall the uniform boundedness of the frequency restriction operator 
$\phi_j(\Lambda_D)$ and some fundamental property of the Besov spaces 
for our purpose of this paper.

\begin{lem}\label{lem:305-1} 
{\rm(}\cite{IMT-2018,Ou_2005,ThOuSi-2002}{\rm)}
Let $1 \leq p \leq \infty$. Then 
\[
\sup _{ j \in \mathbb Z}\| \phi_j(\Lambda_D) \|_{L^p\to L^p}  < \infty . 
\]
In particular, we also have that for any $\varphi \in C_0 ^\infty ((0,\infty))$,   
there exists $C > 0$ such that 
\[
\| \varphi(2^{-j}\Lambda_D )\phi_j(\Lambda_D) f \|_{L^p}
\leq C \| \phi_j(\Lambda_D) f\|_{L^p}
\]
for all $j \in \mathbb Z$ and $f$ with  $\phi_j(\Lambda_D) f \in L^p (\mathbb R^2_+)$. 
\end{lem}

\begin{lem} {\rm(}\cite{IMT-2019}{\rm)} \label{lem:11}
{\rm (i)} {\rm(}Resolution of identity{\rm)} For every $f \in \mathcal X_D '$, we have 
\[
f = \psi(\Lambda _D) f + \sum _{ j =1}^\infty  \phi_j(\Lambda _D ) f  
\quad \text{ in } \mathcal X_D' .
\]
For every $f \in \mathcal Z_D '$, we have 
\[
f = \sum_{ j \in \mathbb Z} \phi_j(\Lambda _D) f \quad 
\text{ in } \mathcal Z_D' .
\]

\noindent 
{\rm (ii)} {\rm(}A characterization of homogeneous spaces as a subspace of 
$\mathcal X_D '${\rm)} 
Let $s < 2/p$ or $(s,p) = (2/p , 1)$. Then 
$\dot B^s_{p,q}(\Lambda)$ is equivalent to 
\[
\Big\{ f \in \mathcal X_D' \, \Big| \, 
  \| f \|_{\dot B^s_{p,q}(\Lambda_D)} < \infty , 
  \quad 
  f = \sum _{ j \in \mathbb Z } \phi_j(\Lambda_D) \text{ in } \mathcal X_D'
\Big\} . 
\]
\end{lem}

Next, we show that 
functions that $\dot B^0_{\infty,1}(\Lambda _D) \cap \dot B^m_{\infty,1}(\Lambda _D)$ 
is included in $C^m (\mathbb R^2_+ \cup \partial \mathbb R^2_+)$ 
and the odd extention for $x_2$ is in $C^m(\mathbb R^2)$. 

\begin{lem}
Let $m = 0,1,2,\cdots$ and $ f \in \dot B^0_{\infty,1}(\Lambda _D) \cap \dot B^m_{\infty,1}(\Lambda _D)$. 
Then $ \partial _x ^\alpha f $ is in $ L^\infty (\mathbb R^2_+) $ and 
is extended to a continuous function in the closure of $\mathbb R^2_+$ 
provided that 
$|\alpha| \leq m$. 
In particular, $f_{odd}$ is regarded as a function belonging to $C^m (\mathbb R^2)$. 
\end{lem}

\noindent 
{\bf Proof. } 
Let $f \in \dot B^0_{\infty,1}(\Lambda _D)$. We write $f$ with the resolution of the identity and consider 
the estimate in $L^\infty (\mathbb R^2_+)$. 
\[
\| f\|_{L^\infty  }
\leq  \sum _{ j \in \mathbb Z} \| \phi_j(\Lambda _D)f \| _{L^\infty} 
= \| f \|_{\dot B^0_{\infty,1}(\Lambda_D)} < \infty. 
\]
Next, we estimate the first derivatives of $\phi_j (\Lambda _D)f$ with each $j \in \mathbb Z$, 
by the derivative estimate for $e^{t\Delta _D}$ and the uniform boundedness 
in Lemma~\ref{lem:305-1} 
\[
\| \nabla  \phi_j(\Lambda _D ) f \|_{L^\infty}
= \| \nabla e^{t_j \Delta _D} (e^{-t_j\Delta _D} \phi_j (\Lambda _D) ) f
\|_{L^\infty}
\leq C 2^{j} \| f \|_{L^\infty} < \infty , 
\quad t_j := 2^{-2j}, 
\]
which implies the uniform continuity of $\phi_j(\Lambda _D ) f$ with respect to $x \in \mathbb R^2_+$. 
The continuity and the following convergence 
\[
\Big\| f - \sum _{ |j| \leq N} \phi_j (\Lambda _D) f \Big\| _{L^\infty}
\leq \Big\| f - \sum _{ |j| \leq N} \phi_j (\Lambda _D) f \Big\| _{\dot B^0_{\infty,1}(\Lambda _D)} 
\to 0 \text{ as } N \to \infty
\]
yield the uniform continuity of $f$, 
and we can then extend $f$ as a continuous function up to the boudary 
$\partial \mathbb R^2_+$. 
We may then have $f \in L^\infty (\mathbb R^2_+) \cap C(\mathbb R^2_+ \cup \partial \mathbb R^2_+)$. 

Let $f \in \dot B^0_{\infty ,1}(\Lambda) \cap \dot B^m_{\infty,1} (\Lambda _D)$ 
for $ m \geq 0 $. Since 
we know $e^{-t_j \Delta_D} \phi_j(\Lambda _D) f \in \dot B^0_{\infty,1}(\Lambda _D)$ 
can be extended to a continuous function on the closure of $\mathbb R^2_+$, 
we write $\phi_j(\Lambda_D)f$ as 
a convolution of the Gauss kernel and an odd function for $x_2$.  
\begin{equation}\label{318-1}
\phi_j(\Lambda _D) f 
= e^{t_j \Delta_D} e^{- t_j \Delta_D} \phi_j(\Lambda _D) f 
= \int_{\mathbb R^2} G_{t_j} (x-y) 
  \Big( e^{-t_j \Delta_D} \phi_j(\Lambda _D) f \Big)_{odd} (y) ~dy ,
\end{equation}
where $t_j = 2^{-2j}$. 
It follows from Lemma~\ref{lem:305-1} that for every mutiindex $\alpha$ with $|\alpha | \leq m$
\[
\begin{split}
\| \partial _x ^\alpha \phi_j(\Lambda _D) f  \|_{L^\infty} 
\leq 
& 
C \| \partial_x^\alpha G_{t_j} \|_{L^1} 
 \| e^{-t_j\Delta_D}\phi_j(\Lambda _D) f \|_{L^\infty}
\\
\leq 
& 
C 2^{|\alpha|j} \| e^{-t_j\Delta_D}\phi_j(\Lambda _D) f \|_{L^\infty}
\\
\leq 
& C2^{|\alpha| j} \| \phi_j(\Lambda _D) f \|_{L^\infty} .
\end{split}
\]
We have from the inequality above that 
\[
\| \partial_x ^\alpha f \|_{L^\infty} 
\leq \sum_{j \in \mathbb Z}
\| \partial _x ^\alpha \phi_j(\Lambda _D) f
\|_{\dot B^0_{\infty,1}(\Lambda_D)} 
\leq C \| f \|_{\dot B^{|\alpha|}_{\infty,1}(\Lambda_D)} ,
\]
$\partial_x^\alpha f$ is extented to a continuous function 
in the closure of $\mathbb R^2_+$ 
and $f_{odd}$ is regarded as a function on $\mathbb R^2$. 
We can also see that $\partial _x ^\alpha f$ is uniformly continuous 
in a similar way to the previous case when $ m= 0$. 

Finally it is easy to see that the odd extention $f_{odd}$ 
is regarded as a function in $C^m(\mathbb R^2)$, since the last right hand side of \eqref{318-1} 
can be regarded as an odd function by taking $x $ in the whole space.  
\hfill $\Box$

\bigskip 

Next lemma reveals that the spectral restriction operator $\varphi(\Lambda_D)$ 
is written by using the Fourier transform and the odd extention, and 
relation between derivatives and smooth odd or even extention. 

\begin{lem}\label{lem:305-2}
{\rm (i)} Let $\varphi \in C_0^\infty ((0,\infty))$ and 
$f \in L^1 (\mathbb R^2_+) + L^\infty (\mathbb R^2_+)$. 
Then 
\begin{equation}\label{305-2}
\begin{split}
\varphi (\Lambda _D) f 
=
&  \int_{\mathbb R^2_+} 
  \left( \mathcal F^{-1}\Big[\varphi (|\xi|)\Big] (x-y) 
       - \mathcal F^{-1} \Big[ \varphi (|\xi|) \Big](x_1-y_1, x_2+y_2) 
  \right) f(y)~dy 
\\
=
& 
 \mathcal F^{-1}\Big[ \varphi(|\xi|) \Big] *  f _{odd} 
 \Big|_{\mathbb R^2_+} .
\end{split}
\end{equation}

\noindent 
{\rm (ii)} Let $f \in \dot B^m_{\infty,1}(\Lambda _D)$ for all $ m = 0,1,2,\cdots$. 
Then 
\[
(\partial _{x_1} f) _{odd} = \partial _{x_1} f_{odd}, \quad 
(\partial_{x_2} f)_{even} = \partial _{x_2} (f_{odd}). 
\]

\end{lem}

\noindent
{\bf Proof. } 
We prove {\rm (i)}. 
Since the support of $\varphi$ is away from the origin, there is 
$\widetilde \varphi \in C_0^{\infty}(0,\infty)$ such that 
\[
\widetilde \varphi(\lambda^2) = \varphi(\lambda ), 
\]
and let us consider $\widetilde \varphi (-\Delta_D)$ instead of $\varphi (\Lambda_D)$, 
whose proof requires essentially same argument. 
Therefore we will show that 
\begin{equation}\label{305-6}
\widetilde \varphi (-\Delta _D ) f = 
\mathcal F^{-1}[\widetilde \varphi (|\xi|^2)] * f_{odd} \Big|_{\mathbb R^2_+} . 
\end{equation}
We can also suppose that 
there exists $j_0 \in \mathbb N$ such that 
\begin{equation}\label{305-3}
f = \sum_{ |j| \leq j_0} \phi_j (\Lambda_D) f, 
\end{equation}
since 
\[
\widetilde \varphi(-\Delta) f 
= \widetilde \varphi(-\Delta)\sum_{|j| \leq j_0} \phi_j(\Lambda _D) f 
\]
for sufficiently large $j_0$. 

We recall the spectral multiplier theorem with the bound by $H^s$ norm 
(see~\cite{Ou_2005,ThOuSi-2002}), and it holds that 
\begin{equation}\label{305-4}
\| \widetilde \varphi (-\Delta _D) \|_{L^\infty \to L^\infty} 
\leq C \| \widetilde \varphi  \|_{H^l (\mathbb R)} 
\text{ for } l > \dfrac{n}{2} .
\end{equation}
Let us fix $l > n/2$. 
Our strategy is to approximate $\widetilde \varphi$ by an polynomial,  
more precisely, we utilize an analytic function. 
For every $\varepsilon > 0$, there exists $\widetilde \varphi_\varepsilon$ such that 
the support of the Fourier transform of $\widetilde \varphi _{\varepsilon} $ 
is compact and 
\[
\| \widetilde \varphi - \widetilde \varphi_{\varepsilon} \|_{H^l (\mathbb R)} < \varepsilon.
\]
By the compactness of the Fourier support, there exists $C>0$ such that 
\[
\| \nabla^\alpha \widetilde \varphi_\varepsilon \|_{L^\infty} 
\leq C ^{|\alpha|} \| \widetilde \varphi_{\varepsilon} \|_{L^\infty},
\]
and the Taylor expansion of $\widetilde \varphi_\varepsilon$ has the 
convergence radius, infinity, for each point. 
We take $\lambda _0 > 0$ and write the Talor expansion 
\[
\widetilde \varphi_\varepsilon (\lambda ) 
= \sum _{k=0}^\infty a_k (\lambda - \lambda _0)^k,
\]
where $a_k$ ($k = 0,1,2,\cdots$) are real numbers 
and the convergence of the series is uniform on each bounded interval. 
It follows by the formula above that 
\begin{equation}\label{305-5}
\begin{split}
& 
\widetilde \varphi(-\Delta_D) f - \mathcal F^{-1}[\widetilde \varphi(|\xi|)] * f_{odd} |_{\mathbb R^2_+}
\\
= 
&  \widetilde \varphi(-\Delta _D) f - \widetilde \varphi_{\varepsilon} (-\Delta _D) f 
+ \sum _{k=0}^\infty a_k (-\Delta _D - \lambda _0)^k  f 
- \mathcal F^{-1}[\widetilde \varphi(|\xi|)] * f_{odd} |_{\mathbb R^2_+} .
\end{split}
\end{equation}
On the first two terms, we have from the boundedness of the spectral multiplier~\eqref{305-4} that 
\[
\| \widetilde \varphi(-\Delta _D) f - \widetilde \varphi_{\varepsilon} (-\Delta _D) f 
\|_{L^\infty}
\leq C \| \widetilde \varphi - \widetilde \varphi_\varepsilon \|_{H^l (\mathbb R)} 
    \| f \|_{L^\infty}
\leq C \varepsilon \| f \|_{L^\infty} . 
\]
On the third and fourth terms, we note that 
\[
-\Delta _D f \in L^\infty 
\quad \text{ if and only if }  \quad 
\lim_{t\to0} \frac{e^{t \Delta _D}f-f}{t} \text{ in } L^\infty (\mathbb R^2_+)  \text{ exists}, 
\]
and we know $-\Delta _D f \in L^\infty (\mathbb R^2_+)$ 
for $f \in L^\infty (\mathbb R^2_+)$ since \eqref{305-3}. 
We also have that 
\[
\lim_{t\to0} \frac{e^{t \Delta _D}f-f}{t} 
=\lim_{t\to0} \frac{e^{t\Delta_{\mathbb R^2}} f_{odd} |_{\mathbb R^2_+}-f}{t} 
= -\Delta f_{odd} |_{\mathbb R^2_+}  
\text{ in } L^\infty (\mathbb R^2_+) 
\text{ exists},
\]
where $\Delta_{\mathbb R^2}$ is the Lalacian on $\mathbb R^2$. 
This implies that 
\[
\sum _{k=0}^\infty a_k (-\Delta _D - \lambda _0)^k  f
=\sum _{k=0}^\infty a_k (-\Delta_{\mathbb R^2}- \lambda _0)^k  f_{odd} 
  \Big| _{\mathbb R^2_+}
= \mathcal F^{-1}[\widetilde \varphi_\varepsilon (|\xi|)] * f_{odd} 
  \Big| _{\mathbb R^2_+} .
\]
This equality and the boudedness of the Fourier multiplier implies that 
\[
\Big\| 
\sum _{k=0}^\infty a_k (-\Delta _D - \lambda _0)^k  f 
- \mathcal F^{-1}[\widetilde \varphi(|\xi|)] * f_{odd} |_{\mathbb R^2_+}
\Big\| _{L^\infty}
\leq C \| \widetilde \varphi_{\varepsilon} - \widetilde \varphi \|_{H^l(\mathbb R)} 
\| f \|_{L^\infty}
\leq C\varepsilon \| f \|_{L^\infty}. 
\]
Since $\varepsilon>0$ is arbitrary, 
we obtain \eqref{305-6} by \eqref{305-5} and the two inequalities above.

The second statement (ii) follows from the equality \eqref{305-2} 
and a symmetric property of the evenness and the oddness with respect to $x_2$ 
taking $x$ in the whole space $\mathbb R^2$. 
In fact, we can write 
\[
\partial _{x_k}f 
= \sum_{j \in \mathbb Z} \partial_{x_k} \phi_j(\Lambda _D) f 
=
\sum_{j \in \mathbb Z} 
\partial_{x_k} 
\mathcal F^{-1}[\phi_j(|\xi|)]* f_{odd} \Big| _{\mathbb R^2_+}, 
\quad k = 1,2,
\]
and notice the radially symmetricity 
$\mathcal F^{-1}[\phi_j(|\xi|)](x) = \mathcal F^{-1}[\phi_j(|\xi|)](|x|) $. 
For $x_1$ derivative, it is easy check that 
\[
\partial_{x_1} \mathcal F^{-1}[\phi_j(|\xi|)]* f_{odd} 
= \int_{\mathbb R^2} \frac{x_1-y_1}{|x-y|} 
  \Big( \partial_{|x-y|} \mathcal F^{-1}[\phi_j(|\xi|)]\Big) (x-y) f_{odd} (y)dy, 
  \quad x \in \mathbb R^2 
\]
is odd for $x_2$, which proves that 
\[
(\partial _{x_1}f) _{odd} 
= \sum_{j \in \mathbb Z} 
\partial_{x_1} 
\mathcal F^{-1}[\phi_j(|\xi|)]* f_{odd}
= \partial _{x_1} f_{odd}  
\quad \text{ in  } \mathbb R^2 .
\]
For $x_2$ derivative, we see that 
\[
\partial_{x_2} \mathcal F^{-1}[\phi_j(|\xi|)]* f_{odd} 
= \int_{\mathbb R^2} \frac{x_2-y_2}{|x-y|} 
  \Big( \partial_{|x-y|} \mathcal F^{-1}[\phi_j(|\xi|)]\Big) (x-y) f_{odd} (y)dy, 
  \quad x \in \mathbb R^2 
\]
is even for $x_2$, and hence, 
\[
(\partial _{x_2}f) _{even} 
= \sum_{j \in \mathbb Z} 
\partial_{x_2} 
\mathcal F^{-1}[\phi_j(|\xi|)]* f_{odd}
= \partial _{x_2} f_{odd}  
\quad \text{ in  } \mathbb R^2 .
\]
\hfill $\Box$

\section{Proof of Theorem~\ref{thm:small}}\label{sec:pf_small}

We prepare two propositions: The first proposition is about the linear estimate 
and the second proposition studies bilinear estimates of $(u \cdot \nabla) \theta$. 
We then prove Theorem~\ref{thm:small}. 

\begin{prop}\label{prop:2} {\rm(}\cite{Iw-2018,Iw-2020}{\rm)}
Let $s \in \mathbb R$, $1 \leq q, q_1 , q_2 \leq \infty$. 

\noindent 
{\rm(i)} {\rm(}Smoothing estimates{\rm)} 
If $s > 0$ and $f \in \dot B^0_{\infty,q} (\Lambda _D)$, then 
\begin{equation}\label{0302-1}
\| e^{-t\Lambda_D} f \|_{\dot B^s_{\infty,q_1}(\Lambda_D)} 
\leq C t^{-s} \| f \|_{\dot B^0_{\infty,q_2}(\Lambda _D)}.
\end{equation}

\noindent 
{\rm(ii)} 
{\rm(}Maximal regularity{\rm)} 
If $f \in \dot B^0_{\infty,q}(\Lambda _D)$, then 
\begin{equation}\label{0301-0}
\begin{split}
\|e^{-t\Lambda} f \|_{B^0_{\infty,q} (\Lambda _D)}
+
\| e^{-t\Lambda} f 
\|_{\widetilde L^1 (0,T ; B^1_{\infty,q}(\Lambda _D)) }
  \leq  C \| f\|_{ B^0_{\infty,q}(\Lambda _D)}. 
\end{split}
\end{equation}
If $f \in L^1 (0,\infty ; \dot B^0_{\infty,1})$, then 
\begin{equation}\label{0301-0(2)}
\Big\| \int_0^t e^{-(t-\tau)\Lambda_D} f(\tau) ~d\tau 
\Big\|_{L^\infty (0,\infty; \dot B^0_{\infty,1}(\Lambda_D)) \cap L^1 (0,\infty ;\dot B^1_{\infty,1}(\Lambda _D)) } 
\leq \| f \|_{L^1 (0,\infty ;\dot B^0_{\infty,1}(\Lambda ))}. 
\end{equation}

\noindent 
{\rm (ii)} 
There exists a positive constant $C$ such that  
\begin{equation}\label{0301-1}
\begin{split}
& 
\| t^{\beta _1 + \beta_2} \partial_t ^{\beta _1} \Lambda _D ^{\beta_2} e^{-t\Lambda_D} f\|_{\dot B^0_{\infty,q} (\Lambda _D)}
+
\| t^{\beta _1 + \beta_2} \partial_t ^{\beta_1} \Lambda _D ^{\beta_2} e^{-t\Lambda_D } f
\|_{\widetilde L^1 (0,\infty ; \dot B^1_{\infty,q}(\Lambda _D)) }
\\
  \leq 
  &
   C^{\beta _1 + \beta_2 } \beta _1! \beta_2 ! 
     \, \| f \|_{\dot B^0_{\infty,q}(\Lambda _D)},
\end{split}
\end{equation}
for all $f \in  \dot B^0_{\infty,q} (\Lambda _D)$ and $\beta _1 , \beta _2 \in \mathbb N \cup \{ 0 \}$. 

\end{prop}

\noindent 
{\bf Remark. } 
The proof of (i) and (ii) can be found in \cite{Iw-2018}. 
The inequality (ii) can be proved in a similar way to the proof of 
Proposition~3.2 in \cite{Iw-2020}. 
Corresponding estimate to \eqref{0301-1} 
has the constant $(\beta_1 + \beta _2)!$, but it should have been 
$\beta _1 ! \beta _2 $ for the proof of the analyticity.  
It is possible to modify the proof to obtain $\beta _1 ! \beta _2 !$, 
by estimating the derivative orders of $t,x$ separatelly.  

\begin{prop}\label{prop:0309}
Let $s \geq 0$. Then for every $f \in \dot B^{0}_{\infty,1}(\Lambda_D) \cap \dot B^s_{\infty,1}(\Lambda_D)$ 
and $g \in \dot B^1_{\infty,1} (\Lambda _D) \cap \dot B^{s+1}_{\infty,1} (\Lambda _D)$ 
\begin{equation}\label{0301-2}
\begin{split}
& 
\Big\| \Big(\nabla ^{\perp}\Lambda _D ^{-1} f \cdot \nabla\Big) g  
\Big\|_{\dot B^s_{\infty,1}(\Lambda _D)} 
\\ 
\leq 
& C ^{s+1}  
\begin{cases}
\| f \|_{\dot B^0_{\infty,1}(\Lambda _D)} \| g \|_{\dot B^1_{\infty,1}(\Lambda _D)} 
& \text{ if } s = 0,
\\
\| f \|_{\dot B^s_{\infty,1}(\Lambda _D)}  \| g \|_{\dot B^1_{\infty,1}(\Lambda _D)}
+ \| f \|_{\dot B^0_{\infty,1}(\Lambda _D)} \| g \|_{\dot B^{s+1}_{\infty,1}(\Lambda _D)} 
& \text{ if } s > 0 .
\end{cases}
\end{split}
\end{equation}

\end{prop}

\noindent 
{\bf Remark. }  
In the whole space case, the inequality \eqref{0301-2} is proved based by Bony paraproduct 
formula~\cite{Bony-1981}, 
\[
\Big(\nabla ^{\perp}(-\Delta _{\mathbb R^2}) ^{-1/2} f \cdot \nabla\Big) g 
= \Big( \sum_{k \geq l+3} + \sum _{l\geq k+3} + \sum _{|k-l| \leq 2} \Big) 
\Big(\nabla ^{\perp}(-\Delta _{\mathbb R^2}) ^{-1/2} 
  f_k \cdot \nabla\Big) g _l,
\]
where $f_k = \mathcal F^{-1}[\phi_j(|\xi|)] *f$ and $g_l = \mathcal F^{-1}[\phi_l(|\xi|)] *g$. 
In fact, the first two terms are estimated by uniformity of the spectral multiplier. 
The third term needs an additional argument by the divergenve free property from $\nabla ^{\perp}$ 
and 
\[
\Big(\nabla ^{\perp}(-\Delta _{\mathbb R^2}) ^{-1/2} f_k \cdot \nabla\Big) g_l
= \nabla \cdot \Big(  \big(\nabla ^{\perp}(-\Delta _{\mathbb R^2}) ^{-1/2} f_k \big) g_l  \Big) .
\]
One can prove the dependence for the constant $C^{s+1}$ in \eqref{0301-2} with respect 
to the regularity number $s \geq 0$ by estimating the paraproduct formula carefully, 
and we left it to the reader. 
We admit the estimate in the whole space case in the proof of Proposition~\ref{prop:0309} below. 

\bigskip

\noindent 
{\bf Proof}.  We may assume that $f,g \in \dot B^{s}_{\infty,1}(\Lambda_D)$ for all $ s \in \mathbb R$, since 
the intersection of $\dot B^{s}_{\infty,1}(\Lambda _D) $ with 
all $s$ is dense in $\dot B^{s}_{\infty,1}(\Lambda _D)$ for 
each $s \in \mathbb R$. 
We write 
\begin{equation}\label{305-7}
\Big(\nabla ^{\perp}\Lambda _D ^{-1} f \cdot \nabla\Big) g  
= - \partial_{x_1} \Big( (\partial_{x_2} \Lambda _D ^{-1} f )g \Big) 
+ \partial_{x_2} \Big( (\partial_{x_1} \Lambda _D ^{-1} f ) g \Big),
\end{equation}
and handle the two terms separately, with Lemma~\ref{lem:305-2} (ii). 
We need the odd extension of the product above to take 
the norm of $\dot B^s_{\infty,1}(\Lambda _D)$. 
For the first term, we write by Lemma~\ref{lem:11}~(i) and $f \in \dot B^s_{\infty,1}(\Lambda _D)$ for all $s \in \mathbb R$
\[
 \Lambda _D ^{-1} f 
 = \sum_{j \in \mathbb Z} \mathcal F^{-1}\left[|\xi|^{-1}\phi_j (|\xi|) \right] * f_{odd} \Big|_{\mathbb R^2_+} 
 = (-\Delta_{\mathbb R^2})^{-1/2} f_{odd} \Big|_{\mathbb R^2_+} 
 \quad \text{in } \mathcal X_D ' ,
 \]
and
\[
\begin{split}
\partial_{x_2} \Lambda _D ^{-1} f 
= 
 \partial_{x_2} (-\Delta_{\mathbb R^2}) ^{-1/2} f_{odd}\Big|_{\mathbb R^2_+} , 
\end{split}
\]
and notice that 
$\partial_{x_2} (-\Delta_{\mathbb R^2}) ^{-1/2} f_{odd}$ is even for $x_2$. 
We then write the odd extention of the above product 
\[
\begin{split}
\bigg( \partial_{x_1} \Big( (\partial_{x_2} \Lambda _D ^{-1} f )g \Big) 
\bigg) _{odd}
=
& 
\partial _{x_1}\bigg( 
\Big( \partial_{x_2} \big((-\Delta_{\mathbb R^2}) ^{-1/2} f_{odd} \big) \Big)
\Big|_{\mathbb R^2_+}  g 
\bigg)_{odd}
\\
=
& 
\partial_{x_1} 
\bigg( \Big( \partial_{x_2} \big((-\Delta_{\mathbb R^2}) ^{-1/2} f_{odd} \big) \Big)
g _{odd}
\bigg) , 
\end{split}
\]
Similarly, we write the second term 
\[
\begin{split}
\bigg( \partial_{x_2} \Big( (\partial_{x_1} \Lambda _D ^{-1} f ) g \Big)
\bigg)_{odd}
=
& 
 \partial_{x_2} \bigg(\Big( (\partial_{x_1} \Lambda _D ^{-1} f ) g \Big)_{even}
\bigg) 
= 
 \partial_{x_2} \bigg( (\partial_{x_1} \Lambda _D ^{-1} f )_{odd} \cdot g _{odd}
\bigg) 
\\
=
& \partial_{x_2} \bigg( (\partial_{x_1} (-\Delta_{\mathbb R^2})^{-1/2} f_{odd} ) \cdot g _{odd}
\bigg) .
\end{split}
\]
It follows from the above two equalities and the bilinear estimates in Besov spaces 
on the whole space $\mathbb R^2$ that 
\[
\begin{split}
\Big\| \Big(\nabla ^{\perp}\Lambda _D ^{-1} f \cdot \nabla\Big) g    
\Big\|_{\dot B^0_{\infty,1}(\Lambda _D)} 
= 
& \sum_{j\in\mathbb Z} 
  \bigg\| \mathcal F^{-1} [\phi_j (|\xi|)  ] * 
        \nabla \cdot 
\bigg( \Big( \nabla ^{\perp} \big((-\Delta_{\mathbb R^2}) ^{-1/2} f_{odd} \big) \Big)
g _{odd}
\bigg) 
  \bigg\|_{L^\infty(\mathbb R^2)}
\\
\leq 
& C \| f_{odd} \|_{\dot B^0_{\infty,1}(\mathbb R^2)} 
    \| g_{odd} \|_{\dot B^1_{\infty,1} (\mathbb R^2)} 
\\
= 
& C \| f \|_{\dot B^0_{\infty,1}(\Lambda_D)} \| g \|_{\dot B^1_{\infty,1}(\Lambda _D)}. 
\end{split}
\]
If $s > 0$, it holds that 
\[
\begin{split}
\Big\| \Big(\nabla ^{\perp}\Lambda _D ^{-1} f \cdot \nabla\Big) g    
\Big\|_{\dot B^s_{\infty,1}(\Lambda _D)} 
\leq C C^s \Big( 
 \| f \|_{\dot B^s_{\infty,1}(\Lambda_D)} \| g \|_{\dot B^1_{\infty,1}(\Lambda_D)} 
 + \| f \|_{\dot B^0_{\infty,1}(\Lambda_D)} \| g \|_{\dot B^{s+1}_{\infty,1}(\Lambda_D)}
 \Big) .
\end{split}
\]
\hfill $\Box$

\bigskip 

\noindent 
{\bf Proof of Theorem~\ref{thm:small}.} 
For the sake of the simplicity, we write 
\[
L^\infty \dot B^0_{\infty,1} := L^\infty (0,\infty ; \dot B^0_{\infty,1}(\Lambda _D)) , \qquad 
L^1 \dot B^1_{\infty,1}  
:= L^1 (0,\infty ; \dot B^1_{\infty,1}(\Lambda _D)).
\]
Let $\Psi$ be the right hand side of the integral equation, 
\[
\Psi (\theta)
:= 
e^{- t \Lambda_D} \theta _0 
  -\int_0^t e^{-(t-\tau)\Lambda_D} 
   \Big( (u \cdot \nabla ) \theta \Big) ~d\tau, \quad 
   u := \nabla^{\perp}\Lambda _D ^{-1}\theta . 
\]

\noindent 
Step 1. (Existence and space analyticity) 
Let a complete metric space $X_\infty$ be defined by 
\[
X_\infty := \{ \theta \in C([0,\infty) , \dot B^0_{\infty,1} (\Lambda _D)) \, | \, 
  \| \theta \|_X \leq \| u_0 \|_{\dot B^0_{\infty,1}(\Lambda _D) }\} , 
\]
where 
\[
\| \theta \|_X := 
\sup _{\beta = 0,1,2,\cdots} 
\frac{
\| t^{\beta} \Lambda _D ^{\beta} \theta 
\|_{L^\infty\dot B^0_{\infty,1}\cap L^1\dot B^1_{\infty,1}} 
}
{C_0 ^{2\beta +1} \beta !}, 
\]
for some large constant $C_0$, 
with the metric 
\[
d(\theta, \tilde \theta):= \| \theta-\tilde \theta\|_{L^\infty (0,\infty ; \dot B^0_{\infty,1}(\Lambda _D)) 
 \cap L^1 (0,\infty ; \dot B^1_{\infty,1}(\Lambda _D))} . 
\]
The main point is to prove that 
\begin{equation}\label{0303-1}
\| \Psi(\theta) \|_X \leq C \| \theta_0 \|_{\dot B^0_{\infty,1}(\Lambda_D)} + C \| \theta  \|_X^2  , 
\text{ for } \theta \in X_\infty
\end{equation}
which implies $\Psi (\theta) \in X_\infty$ provided that $u_0$ is small in $\dot B^0_{\infty,1}(\Lambda_D)$. 
Hereafter, we estimate, supposing the smallness. 

By \eqref{0301-1}, the linear part is estimated by 
\[
\| e^{-t\Lambda _D} \theta _0 \|_{X}
\leq C 
\left( 
\sup _{\beta= 0,1,2,\cdots} 
\frac{ C^{\beta} \beta !}{C_0 ^{2\beta + 1} \beta!}
\right) 
\| \theta_0 \| _{\dot B^0_{\infty,1}(\Lambda_D)},
\]
and we have the finiteness of the supremum 
in the right hand side for large $C_0$.  
We turn to estimate the nonlinear term, dividing the interval in half. 

For  the first-half time integral, 
it follows that by the smoothing effect \eqref{0302-1} 
\[
\begin{split}
& 
\Big\| t^{\beta _2} \Lambda _D^{\beta_2 }
  \int_0^{t/2} e^{-(t-\tau)\Lambda_D} \Big( (u \cdot \nabla ) \theta  \Big)  ~ d\tau 
\Big\|_{\dot B^0_{\infty,1}(\Lambda_D)}
\\
\leq 
&C ^{\beta _2} \beta _2 ! \cdot  t^{\beta_2} \int_0^{t/2} (t-\tau)^{-\beta_1} 
  \| (u\cdot \nabla) \theta \|_{\dot B^0_{\infty,1}} d\tau 
\\
\leq 
& C^{\beta_2} \beta _2 !\int_0^{t/2} \| (u\cdot \nabla) \theta \|_{\dot B^0_{\infty,1}} d\tau ,
\end{split}
\]
and by the smoothing effect \eqref{0302-1} and maximal regularity \eqref{0301-0(2)} 
\[
\begin{split}
& 
\Big\| t^{\beta _2} \Lambda _D^{\beta_2 }
  \int_0^{t/2} e^{-(t-\tau)\Lambda_D} \Big( (u \cdot \nabla ) \theta  \Big)  ~ d\tau 
\Big\|_{L^1(0,\infty ;\dot B^1_{\infty,1}(\Lambda_D))}
\\
\leq 
&C ^{\beta _2} \beta _2 ! 
\Big\| t^{\beta_2} \int_0^{t/2} \Big(\dfrac{t-\tau}{2}\Big)^{-\beta_1} 
  \| e^{-\frac{t-\tau}{2}\Lambda_D} (u\cdot \nabla) \theta \|_{\dot B^1_{\infty,1}} d\tau 
\Big\|_{L^1 (0,\infty)} 
\\
\leq 
& C^{\beta_2} \beta _2 ! \int_0^{\infty} \| (u\cdot \nabla) \theta \|_{\dot B^0_{\infty,1}} d\tau . 
\end{split}
\]
These inequalities above and the bilinear estimate \eqref{0301-2} imply that 
\begin{equation}\notag 
\begin{split}
& 
\Big\| t^{\beta _2} \Lambda _D^{\beta_2 }
  \int_0^{t/2} e^{-(t-\tau)\Lambda_D} \Big( (u \cdot \nabla ) \theta  \Big)  ~ d\tau 
\Big\|_{L^\infty \dot B^0_{\infty,1} \cap L^1 \dot B^1_{\infty,1} }
\\
\leq 
& C^{\beta_2} \beta _2 ! \int_0^\infty 
 \| \theta \|_{\dot B^0_{\infty,1}(\Lambda_D)} \| \theta \|_{\dot B^1_{\infty,1}(\Lambda _D)} d\tau 
 \leq C ^{\beta _2 } \beta _2 ! \| \theta \|_{X}^2 . 
\end{split}
\end{equation}
As for the second-half time integral with $\beta_1 = 0$, by maximal regularity \eqref{0301-0(2)} and the bilinear estimate \eqref{0301-2} 
give that 
\begin{equation}\notag 
\begin{split}
& 
\Big\| t^{\beta _2} \Lambda _D^{\beta_2 }
  \int_{t/2}^t e^{-(t-\tau)\Lambda_D} \Big( (u \cdot \nabla ) \theta  \Big)  ~ d\tau 
\Big\|_{L^\infty \dot B^0_{\infty,1} \cap L^1 \dot B^1_{\infty,1} }
\\
\leq 
&C\int_{0}^\infty  (2s)^{\beta _2} \Big\| \Lambda _D^{\beta_2 }  \Big( (u \cdot \nabla ) \theta  \Big) \Big\|_{\dot B^0_{\infty,1}(\Lambda_D)}  ~ d\tau 
\\
\leq 
&C^{\beta _2}\int_{0}^\infty s  ^{\beta_2}  
\Big( \| \theta \|_{\dot B^{\beta_2}_{\infty,1}(\Lambda_D)} \| \theta\|_{\dot B^1_{\infty,1}(\Lambda_D)} 
   +  \| \theta \| _{\dot B^0_{\infty,1}(\Lambda_D)} \| \theta \|_{\dot B^{\beta_2 + 1}_{\infty,1}(\Lambda_D)} 
\Big)~d\tau 
\\
\leq 
& C ^{\beta_2}\beta_2 ! \cdot \| \theta \|_X^2
\end{split}
\end{equation}
Therefore, we obtain \eqref{0303-1} for large $C_0$, 
and we can also prove that for $\theta , \tilde \theta \in X_\infty$ 
\begin{equation}\label{0304-2}
d(\theta, \tilde \theta)
\leq C (\| \theta \| _X + \| \tilde \theta \|_X ) d(\theta, \tilde \theta)
\leq \frac{1}{2} d(\theta, \tilde \theta),
\end{equation}
provided that $\| u_0 \|_{\dot B^0_{\infty,1}(\Lambda_D)}$ is sufficiently small. 
The fixed point argument yields the existence of the solution, and 
the analyticity for $x$ holds by the estimate of derivatives by factorials. 


\bigskip 

\noindent 
Step 2. (With time analyticity) By the result of Step 1, 
we have  $\Lambda ^\beta \theta (t) \in \dot B^0_{\infty,1}(\Lambda_D)$ 
for all $\beta = 1,2, \cdots $ if $t > 0$. 
Then we know the solution becomes smooth, and let us  assume the initial data $\theta_0 \in \dot B^s_{\infty,1}(\Lambda_D)$ 
for all $s \geq 0$. Our argument is to construct a solution for smooth data.

Let a complete metric space $Y_\infty$ be defined by 
\[
Y_\infty := \{ \theta \in C([0,\infty) , \dot B^0_{\infty,1} (\Lambda _D)) \, | \, 
  \| \theta \|_Y \leq \| u_0 \|_{\dot B^0_{\infty,1}(\Lambda _D) }\} , 
\]
where 
\[
\| \theta \|_Y := 
\sup _{\beta_1 , \beta _2 =0,1,2,\cdots} 
\frac{ (1+\beta_1 + \beta _2)^4 
\|t^{\beta _1 + \beta _2} \partial _t^{\beta _1} \Lambda _D ^{\beta_2} \theta 
\|_{L^\infty\dot B^0_{\infty,1}\cap L^1\dot B^1_{\infty,1}} 
}
{C_0 ^{3\beta_1 + 2\beta _2 -1} \beta _1! \beta _2!}, 
\]
for some large constant $C_0$, 
with the metric 
\[
d(\theta, \tilde \theta):= \| \theta-\tilde \theta \|_{L^\infty (0,\infty ; \dot B^0_{\infty,1}(\Lambda _D)) 
 \cap L^1 (0,\infty ; \dot B^1_{\infty,1}(\Lambda _D))} . 
\]
Let us focus on the following estimate. 
\[
\| \Psi(\theta) \|_Y \leq C \| \theta_0 \|_{\dot B^0_{\infty,1}(\Lambda_D)} + C \| \theta  \|_Y^2  , 
\text{ for } \theta \in Y_\infty
\]
To this end, we will apply an induction argument for time derivatives. 
For $\beta = 0,1,2,\cdots $ we write 
\[
\| \theta \|_{Y _{\leq \beta}}:= 
\sup _{0 \leq \beta_1 \leq \beta , \, \beta _2 = 0,1,2,\cdots} 
\frac{ (1+\beta_1 + \beta _2)^4 
\|t^{\beta _1 + \beta _2} \partial _t^{\beta _1} \Lambda _D ^{\beta_2} \theta 
\|_{L^\infty\dot B^0_{\infty,1}\cap L^1\dot B^1_{\infty,1}} 
}
{C_0 ^{3\beta_1 + 2\beta _2 +1} \beta _1! \beta _2!}. 
\]

When $\beta = 0$, the estimate is essentially proved in Step 1, modifying 
$C_0$ larger to lead to a boundedness with an additional factor $(1+\beta_1+\beta_2)^4$. 

Let $\beta\geq 1$ and assume that 
\[
\| \Phi(\theta) \|_{Y_{\leq \beta-1}} \leq \| u_0 \|_{\dot B^0_{\infty,1}} . 
\]
The equation for $t^{\beta} \partial _t^{\beta} \Psi(\theta)$ is 
\[
\partial _t \big( t \partial_t^{\beta} \Psi(\theta) \big) 
+ \Lambda _D (t^{\beta} \partial _t^{\beta} \Psi(\theta)) 
- \beta t^{\beta-1}\partial_t ^{\beta}\Psi (\theta) 
+ t^{\beta} \partial _t^{\beta} (u\cdot \nabla)\theta  = 0 ,
\]
the initial data of $t^{\beta} \partial _t^{\beta} \Psi(\theta)$ 
is zero because of our smooth $\theta _0$, and we have the integral equation 
\[
\begin{split}
t^{\beta} \partial _t^{\beta} \Psi(\theta)(t) 
= \int_0^t e^{-(t-\tau)\Lambda_D} 
\Big( \beta \tau ^{\beta-1}\partial_\tau^{\beta} \Phi(\theta)(\tau) 
  - \tau^{\beta} \partial_\tau^{\beta} \Big( (u \cdot\nabla)\theta \Big)\Big) d\tau . 
\end{split}
\]
Maximum regularity \eqref{0301-0(2)} implies that 
\begin{equation}\label{0304-1}
\| t^{\beta} \partial _t^{\beta} \Psi(\theta)(t)\| _{L^\infty\dot B^0_{\infty,1} \cap L^1 \dot B^1_{\infty,1}}
\leq C \int_0^\infty  
\Big( \beta \tau ^{\beta-1} \| \partial _\tau ^\beta  \Phi(\theta) \|_{\dot B^0_{\infty,1}} 
    + \tau ^{\beta} \| \partial_\tau ^{\beta} (u \cdot \nabla) \theta\|_{\dot B^0_{\infty,1}} \Big) d\tau . 
\end{equation}
For the first integrand, we write the equation 
$\partial _\tau \Phi(\theta) = - \Lambda _D \Phi(\theta) - (u \cdot \nabla) \theta$ 
and apply the assumption of the indection to the first term in the right hand side, 
and the Leibniz rule, the bilinear esimate \eqref{0301-2} 
for the second term. We then have that 
\[
\begin{split}
& 
C\int_0^\infty  
 \beta \tau ^{\beta-1} \| \partial _\tau ^\beta  \Phi(\theta) \|_{\dot B^0_{\infty,1}} d\tau 
\\
 \leq 
& C \beta \int_0^t \tau ^{\beta -1} \| \partial _{\tau}^{\beta -1} \Psi(\theta) \|_{\dot B^1_{\infty,1}} d\tau 
 + C \beta \int_0^\infty \tau ^{\beta-1}
   \sum_{\gamma = 0}^{\beta -1} \frac{(\beta-1) !}{(\beta-1  -\gamma)! \gamma !} 
    \| \partial_\tau ^{\beta-1 - \gamma}\theta \|_{\dot B^0_{\infty,1}} \| \partial _\tau ^\gamma \theta \|_{\dot B^0_{\infty,1}} d\tau 
\\
 \leq 
& C \beta \cdot \frac{C_0^{3 (\beta -1)+2+1} (\beta _1 -1)!}{(1+\beta)^4} \| u_0 \|_{\dot B^0_{\infty,1}} 
 + C \sum_{\gamma = 0}^\beta \frac{C_0^{3 (\beta-1) + 2} \beta !}{(1+ \beta -1- \gamma)^4 (1+\gamma )^4} 
    \| u_0 \|_{\dot B^0_{\infty,1}}^2 
\\
\leq 
& C \Big( \frac{1}{C_0} + \| u_0 \|_{\dot B ^0_{\infty,1}} \Big) \frac{C^{3\beta +1} \beta !}{(1+\beta )^4} \| u_0 \|_{\dot B^0_{\infty,1}}. 
\end{split}
\]
As for the second term of the right hand side of \eqref{0304-1}, a similar argument to the second estimate above implies that 
\[
\begin{split}
C \int_0^\infty  
   \tau ^{\beta} \| \partial_\tau ^{\beta} (u \cdot \nabla) \theta\|_{\dot B^0_{\infty,1}} d\tau 
\leq C \frac{C_0 ^{3\beta +1} \beta !}{(1+\beta)^4} \cdot C _0\| u_0 \|_{\dot B^0_{\infty,1}}^2 .
\end{split}
\]
Therefore we obtain that 
\[
\| \Phi(\theta) \|_{Y_{\leq \beta}} \leq \| u_0 \|_{\dot B^0_{\infty,1}} . 
\]
This together with \eqref{0304-2} allows us to apply the fixed point argument and we have the solution 
analytic in space and time. 

We can also prove the uniqueness analogously to the paper~\cite{Iw-2015} 
by introducing odd extention with respect to $x_2$, where the uniqueness in 
$C([0,\infty),\dot B^0_{\infty,1}(\Lambda_D)) \cap L^1(0,\infty ; \dot B^1_{\infty,1}(\Lambda _D))$ is proved without smallness. 
It follows from the uniqueness that the solutions in Step 1 and Step 2 concide. 
Overall, we complete the proof of Theorem~\ref{thm:small}. 
\hfill $\Box$ 

%

\section{Proof outline of Theorem~\ref{thm:1}}  

In this section, we explain proof outline, since the argument is similar to 
the whole space case, once we introduce odd extention of the equations. 
The essence of how to introduce the odd extetnion is found in section 2. 
We write the odd extention of the equations. 
\[
\partial_t \theta_{odd} + (-\Delta _D)^{1/2}  \theta_{odd} 
+ \Big( (u \cdot \nabla) \theta \Big)_{odd}  = 0. 
\]
The odd extention of the nonlinear term becomes 
\[
\Big( (u \cdot \nabla) \theta \Big)_{odd}
= (\nabla^\perp (-\Delta_{\mathbb R^2})^{-1/2} \theta_{odd} \cdot \nabla ) \theta _{odd}.
\]
Therefore, we have the equation in the whole space, and it is possible to 
apply the argument in the paper~\cite{Iw-2020,WZ-2011} with a certain modification for the 
low spectral component due to the second condition of \eqref{317-1} for the boundedness of 
the Riesz transform. 
We then obtain a local-in-time unique solution in the whole space, 
analitic in space time for some $T_0 >0$. The restriction to the half space 
gives the solution such that 
\[
\theta \in C([0,T_0] , B^0_{\infty,\infty}(\Lambda_D)) \cap 
 \tilde L^1 (0,T_0 ; \dot B^1_{\infty,\infty }(\Lambda _D))  ,
\]
\[
\Big( 1 - \sum_{j \geq 0} \phi_j (\Lambda _D) \Big)\theta \in C((0,T_0], \dot B^0_{\infty,1} (\Lambda_D)) , 
\]
and the uniqueness of the solution in the half space is also proved 
by odd extention and the argument for the whole space.

We turn to prove the global-in-time regularity, based on the nonlinear maximum principle 
(see~\cite{CV-2012} and also~\cite{Con-2017}). 
Let $\theta$ be the solution constructed as above. 
We notice that for almost every $t > 0$, 
$\theta (t) \in B^a_{\infty,1}(\Lambda_D)$, $0 < \alpha < 1$, 
and $T_0 >0$ for the local existence is taken such that 
\[
\| \theta \|_{\tilde L^1 (0,T_0 ; \dot B^1_{\infty,\infty}(\Lambda_D))}  \ll 1,
\]
and that we can extend the existence time as far as such kind of smallness holds 
for the linear solution. 
We fix a time $t_0 \in (0,T_0]$ and consider the data $\theta (t_0)$. 
We here utilize the uniform bounds in H\"older spaces (see Theorem~3.1 in~\cite{Con-2017}), 
and this implies that 
if $0< a \ll 1/ \| \theta(t_0) \|_{L^\infty}$, then  
\begin{equation}\label{318-2}
\sup_{x\not = y} 
\frac{|\theta(t,x) - \theta(t,y)|}{|x-y|^a} 
\leq \sup_{x \not=y} \frac{|\theta(t_0,x) - \theta (t_0,y)|}{|x-y|^a}
=: M_{t_0} , 
\quad \text{ for all }t \geq t_0 ,
\end{equation}
as far as the solution exists and is smooth. 
If we consider $\theta(t_0)$ as a data, 
we have on a short time interval $[0,\delta]$ that 
\[
\| e^{-t\Lambda_D} \theta (\tilde t \,) \|_{\tilde L^1 (0,\delta ; \dot B^1_{\infty,\infty} (\Lambda_D))}
\leq \| e^{-t\Lambda_D} \theta (\tilde t \,) \|_{L^1 (0,\delta ; \dot B^1_{\infty,1} (\Lambda_D))}
\leq C \delta^{a} \| \theta(\tilde t \,) \|_{\dot B^a_{\infty,\infty}(\Lambda_D)} 
\leq C \delta^a M_{t_0}.
\]
Let us take $\delta $ sufficientlly small, and we see that 
it is possible to extend the existence time of the solution longer than $[0,t_0]$, 
which is $[0, t_0 + \delta]$. Since the bound~\eqref{318-2} is independent 
of $t$, it is possible to repeat this procedure and we have the existence time 
$[0,t_0+n\delta]$ for all $n = 1,2, \cdots$,  
which proves the global-in-time regularity of the solution. 
\hfill$\Box$

\vskip15mm

\noindent
{\bf Acknowledgements. }
The author was supported by the Grant-in-Aid for Young Scientists (A) (No.~17H04824)
from JSPS.

\vskip3mm 
%
\noindent 
{\bf Conflict of Interest. }
The author declares that he has no conflict of interest. 
%
%

\end{document}